\documentclass[11pt,twoside]{article}
\usepackage{times}
\usepackage{amsmath}
\usepackage{amsfonts}
\usepackage{amssymb}
\usepackage{amsthm} 

\newtheorem{theorem}{Theorem}[section]

\newtheorem{proposition}[theorem]{Proposition}

\theoremstyle{definition}

\theoremstyle{remark}

\allowdisplaybreaks
\newcommand{\head}[1]{{\noindent \small Proceedings of the International Workshop\\
Future Directions in Difference Equations.\\
June 13-17, 2011, Vigo, Spain.\\
{\sc Pages} #1}
\vspace{15pt}}
\newcommand{\heading}[4]{
\begin{center}
{\Large\bfseries\boldmath #1}\\[\baselineskip]
{\Large #2}\\[0.5\baselineskip]
{\normalsize #3}
\\[0.5\baselineskip] {\small\ttfamily #4}
\end{center}
}
\newcommand{\headline}[2]{\pagestyle{myheadings}\markboth{\centerline{\sc #1}} {\centerline{\sc #2}}}

\begin{document}

\head{ip--fp}


\heading
{The Dynamical Degrees of a Mapping}
{Eric Bedford}
{Indiana University}
{bedford@indiana.edu}
\thispagestyle{plain}


\headline{Eric Bedford}{Dynamical Degrees}
\setcounter{section}{0}
\setcounter{equation}{0}

\begin{abstract}
Let $f:X\to X$ be a rational mapping in higher dimension.  The complexity of $(f,X)$ as a dynamical system is measured by the dynamical degrees $\delta_p(f)$, $1\le p\le {\rm dim}(X)$.  We give the definition of the dynamical degrees show how they are computed in certain cases.  For instance, we show that if the dynamical degree of an automorphism of a K\"ahler manifold is greater than one, then it must be irrational.
\end{abstract}

\section{Dynamical degree}
Let us start by discussing  automorphisms of ${\bf C}^2$. 
We say that   
$$f(x,y) = (f_1(x,y),f_2(x,y)):{\bf C}^2\to{\bf C}^2$$  
is a polynomial mapping if the coordinate functions $f_1$ and $f_2$ are polynomials, and we define the degree of $f$ as  ${\rm deg}(f):=\max({\rm deg}(f_1),{\rm deg}(f_2))$.  The degree is not invariant under conjugation.   That is, if $L$ is linear, then the ${\rm deg}(L)=1$, but if $f$ is a polynomial automorphism, then in general ${\rm deg}(f\circ L\circ f^{-1})\ge 1$, and with suitable choice of $f$, this degree can be arbitrarily large.  The behavior of $deg$ under composition is ${\rm deg}(f\circ g)\le {\rm deg}(f){\rm deg}(g)$.  Thus we may define the dynamical degree as
$$\delta(f):=\lim_{n\to\infty} {\rm deg}(f^n)^{1/n}$$
It follows that $\delta(f)=\delta(h^{-1}\circ f\circ h)$, so the dynamical degree is invariant under conjugation.  The condition $\delta>1$ corresponds to exponential growth of degree under iteration, and this may be viewed as ``degree complexity.''   Let us consider two examples:
$$h(x,y)=(y,\varphi(y)-\alpha x),\ \ \ \ k(x,y)=(x,y+\varphi(x))  \eqno(1.1)$$
where $\varphi$ is a monic polynomial.  We see that the iterative behavior of the two maps in (1.1) is rather different: $\delta(h)={\rm deg}(\varphi)$, and $\delta(k)=1$.   The following result from \cite{[FM]} gives a satisfying characterization of the situation for polynomial automorphisms of ${\bf C}^2$:
\begin{theorem} If $f$ is a polynomial automorphism of ${\bf C}^2$ with $\delta(f)>1$, then $f$ is conjugate to a map of the form $h_1\circ\cdots\circ h_j$, where $h_i=(y,\varphi_i(y)-\alpha_ix)$.  In particular, $\delta(f)={\rm deg}(\varphi_1)\cdots{\rm deg}(\varphi_j)$ is an integer.  
\end{theorem}

\noindent The maps $h_i$ that appear in the Theorem are called generalized H\'enon maps. The H\'enon representation achieves minimal degree, and this representation is an essentially unique representative of the conjugacy class.  Thus if we have a H\'enon representative, we know the dynamical degree.  As will be seen in Theorem 6.1 below, the fact that $\delta(f)$ is an integer prevents $f$ from being conjugate to a compact surface automorphism. 

Now let us consider maps of projective space.  Let $(f_0,\dots,f_k)$ be a  $k+1$-tuple of polynomials which are homogeneous of degree $d$.  We may assume that the $f_i$ have no common factor.  The set ${\cal I}(f):=\{x\in {\bf P}^k:f_0(x)=\cdots = f_k(x)=0\}$ (which is possibly empty) has codimension at least 2.  Then $f=[f_0:\cdots : f_k]:{\bf P}^k-{\cal I}(f)\to {\bf P}^k$ is holomorphic.  At each point  $p\in {\cal I}(f)$, however, $f$ is discontinuous and in fact ``blows up''  $p$ to a set of positive dimension.   A topological fact is that the cohomology groups $H^2({\bf P}^k;{\bf Z})$ and $H^{1,1}({\bf P}^k;{\bf Z})$ are both isomorphic to the Picard group $Pic(X)$.  The Picard group is the set $Div(X)/\sim$ of integral divisors modulo linear equivalence.  That is, a divisor $D$ is linearly equivalent to zero if $D=div(h)$, where $h$ denotes a rational (or meromorphic) function $h$ on $X$, and $div(h)=Zeros(h)-Poles(h)$ is the associated divisor.  $Pic({\bf P}^k)$  is generated by the class of a hyperplane $H=\{\sum c_j x_j=0\}$.   To see this, suppose that $V=\{P=0\}$ is the zero set of a polynomial of degree $m$, then for $0\le j\le k$, $h:=P/x_j^m$ is a well defined rational function, which shows that $[V]=m [H]$ in $Pic$.  The action of $f^*$ on $Pic$ is composition: $f^*\{P=0\}=\{P\circ f=0\}$, so  $f^* [H]=d\cdot [H]$.  

More generally, if $\pi:X\to{\bf P}^k$ is a blowup space, then we have the induced map  $f_X:=\pi^{-1}\circ f\circ \pi$ on $X$.  We have well-defined pullback maps $f^*$ on $H^{1,1}({\bf P}^2)$ and $f_X^*$ on $H^{1,1}(X)$.   We can use $f^*$  to define the degree of $f$.  We can use either $f^*$ or $f_X^*$ to define the dynamical degree:
$$\delta(f)=\lim_{n\to\infty}|| (f^n)^*||^{1/n}    \eqno(1.2)$$
where $||\cdot ||$ denotes any norm on $H^{1,1}(X)$, $H^2(X)$, or in nice cases, $Pic(X)$.  

In particular if $X$ is a compact manifold,  (1.2) can be used to define $\delta(f)$ for any meromorphic map $f:X\to X$.   The following is evident:
\begin{proposition}  If  $(f^n)^*=(f^*)^n$ on $H^{1,1}$ for $n>0$, then $\delta(f)$ is the spectral radius of $f^*$, i.e.,  the modulus of the largest eigenvalue of $f^*$.   In this case, $\delta(f)$ is an algebraic integer.
\end{proposition}

\section{ Finding automorphisms by blowing up space}
Let us illustrate this with maps of the form
$$f_{a,b}(x,y)=\left (y, {y+a\over x +b}\right)$$
for fixed constants $a$ and $b$.  This family is conjugate (via affine transformations) to the family $F_{\alpha,\beta}(x,y)=\left(y,{y/ x}\right) + (\alpha,\beta)$, and we are free to work with the maps in either form.   $f_{a,b}$ is a birational map of the plane, and we may extend $f_{a,b}$ to a compactification of the plane.  We start by extending it to the projective space ${\bf P}^2 = \{[x_0:x_1:x_2]\}$ with $(x,y)\leftrightarrow [1:x:y]$. Thus ${\bf P}^2 = {\bf C}^2 \cup L_\infty$, where $L_\infty=\{x_0=0\}$ is the line at infinity.  In homogeneous coordinates we have
$$f_{a,b}[x_0:x_1:x_2]= [x_0(x_1+b x_0):x_2(x_1+b x_0): x_0(x_2+a x_0)].$$

In order to understand the map $f_{a,b}$, we will try to see whether there is a ``better'' compactification.  We start by observing that there is a triangle of lines which are mapped to points:
$$ L_\infty=\{x_0=0\}    \to   e_1:=[0:1:0],  \ \    \{x+b=0\} = \{bx_0+x_1=0\}\to e_2:=[0:0:1],$$
$$  \{y+a=0\}=\{ax_0+x_2=0\}\to q:= (-a,0) = [1:-a:0] $$
We have given the lines of the triangle both in coordinates $(x,y)$ on ${\bf C}^2$ and $[x_0:x_1:x_2]$ on ${\bf P}^2$.
The points $e_1$, $e_2$ and $p:=(-b,-a)$ are indeterminate.  The point $e_2$, for instance, is contained in both $\{x+b=0\}$ and $L_\infty$, so it must blow up to a connected set containing the images of both of these lines.  In this case we have the simplest possibility:  $e_2$ blows up to $\{x_0=0\}$, the line through $e_2$ and $e_1$.  

We describe the operation of blowing up the origin $(0,0)\in {\bf C}^2$.  We define 
$$\widehat{ {\bf C}^2} = \{(x,\xi)= ((x_1,x_2),[\xi_1:\xi_2])\in {\bf C}^2\times{\bf P}^1: x_1\xi_2=x_2\xi_1\}$$ 
and $\pi(x,\xi)=x$.  We say that $\pi:\widehat{{\bf C}^2}\to{\bf C}^2$ is the blowup map, and the blowup space  $\widehat{ {\bf C}^2} $ is a (smooth) complex manifold with the properties: $E:=\pi^{-1}(0,0)$ is equivalent to ${\bf P}^1$, and $\pi:\widehat{ {\bf C}^2} -E\to {\bf C}^2-(0,0)$ is biholomorphic.  $\widehat{ {\bf C}^2} $ is covered by the open sets $\{\xi_j\ne0\}$, $j=1,2$.  If $\xi_1\ne0$, then we may suppose that $\xi_1=1$ and represent this open set by the coordinate chart $ {\bf C}^2\ni (t,\eta)\to (x,\xi)$, where $x=(t,t\eta)$ and $\xi=[1:\eta]$.   In this coordinate chart, we have $E\cap\{\xi_1\ne0\}=\{t=0\}$.

The blowup is a local operation, and we may construct a manifold $\pi:X\to{\bf P}^2$ by blowing up ${\bf P}^2$ at the points $e_1$ and $e_2$.  
 Here we use the notation $E_j=\pi^{-1}e_j$.   The blowup space $X$ is defined by the properties:  $\pi:X-(E_1\cup E_2)\to {\bf P}^2-\{e_1,e_2\}$ is biholomorphic, and $E_j\cong {\bf P}^1$, for $j=1,2$.  To work in a coordinate chart at $E_2$ we let $\tilde\pi:X\to {\bf P}^2$ be given by $\tilde \pi((x_0,x_1),[\xi_0:\xi_1]) = [x_0:x_1:1]$ be the blowup map over $(x_0,x_1)=(0,0) = [0:0:1]$.  The coordinate chart for $\xi_0\ne0$ is given by  ${\bf C}^2\ni (t,\eta)\to (x,\xi)$ with $x=[t:t\eta:1]$.  Thus the inverse is given by $\tilde\pi^{-1}[x_0:x_1:1] = (t=x_0, \eta=x_1/x_0)$.
 
 Since $\pi$ is a birational map, we have an induced map $f_X:=\pi^{-1}\circ f\circ\pi:X\to X$.  Now we show that the map $f_X$ sends $\{x+b=0\}$ to $E_2$.  For this we write
 $$f:{\bf C}^2\to {\bf P}^2, \ \ f(x,y) = \left[ 1:y:{y+a\over x+b}  \right ] =   \left [{x+b\over y+a}:{y(x+b)\over y+a}:1 \right ]$$
 so $\tilde\pi^{-1}f(x,y) = (t = {(x+b)/( y+a )}, \eta= y)$.  This means that $\{x+b=0\}$ is taken to $\{t=0\}$, i.e., to $E_2$.
 
 A similar computation shows that $f_X$ is a smooth mapping from $E_2$ to $L_\infty=\{x_0=0\}$.  This time we write $\tilde\pi(t,\eta) = [t:t\eta:1] = [1:\eta:t^{-1}]$.  Thus  we have
 $$f_X: \ \ (t,\eta)\mapsto f(\tilde\pi (t,\eta)) = f(\eta,t^{-1}) = \left [1:t^{-1}: {t^{-1}+ a\over \eta+b}\right ] = \left [t:1: {1+at\over \eta+b} \right ]$$
 Thus $f_X$ takes $E_2=\{t=0\}$  to $\{x_0=0\}$, and $f_X$ is smooth for $\eta\ne-b$.
 
 If $p\in{\bf P}^2-\{e_1,e_2\}$, we write $p$ for its image $\pi^{-1}p$ in $X$ and we let $\{y+a=0\}$ denote the closure in $X$ of the  image $\pi^{-1}\{y+a=0\}$.  Arguing as above, we find that
 $\{x+b=0\}\to E_2\to L_\infty\to E_1$, and:
 
\begin{proposition} The only indeterminate point for $f_X$ is $p$, and the only exceptional curve (i.e., the only curve which maps to a point) is $\{y+a=0\}$.
\end{proposition}

Now we define a subset of parameter space
$${\cal V}_n:=\{(a,b)\in{\bf C}^2: f_X^n(q)=p\}=\{(a,b)\in{\bf C}^2: f_{a,b}^n(-a,0)=(-b,-a)\}$$
The following is from \cite{[BK2]}:
\begin{theorem}  Fix $n\ge 0$.  Then $(a,b)\in {\cal V}_n$ if and only if there is a space $\pi:Y\to X$ such that $f_Y$ is an automorphism of $Y$. 
\end{theorem}

Suppose that $(a,b)\in{\cal V}_n$.  Define $Q_j:=f_X^j(q)$ for $0\le j\le n$.  Now let $\pi:Y\to X$ denote the manifold obtained by blowing up the points $q_0,q_1,\dots,q_n$.  We write $Q_j:=\pi^{-1}q_j$.  If we write local charts as we did for the case $\{x+b=0\}$, we see that the set $\{y+a=0\}$ is not exceptional for $f_Y$.  Similarly, working as we did at $E_2$ above, we see that $f_Y$ is not indeterminate at $P=Q_n$.  We saw already that $f_X$ is a local diffeomorphism at all the intermediate points $q_j$, so $f_Y$ is a local diffeomorphism at $Q_j$.

\section{Finding the degree}
If $X$ is a space obtained by blowing up ${\bf P}^2$, then the cohomology groups $H^2(X;{\bf Z})$ and $H^{1,1}(X;{\bf Z}):=H^{1,1}(X;{\bf C})\cap H^2(X;{\bf Z})$ are both isomorphic to the Picard group $Pic(X)$.  The Picard group is the set $Div(X)/\sim$ of integral divisors modulo linear equivalence.    It is a standard fact that if $\pi:X\to{\bf P}^2$ is the blow up of ${\bf P}^2$ at distinct points $p_1,\dots,p_N$, then a ${\bf Z}$-basis for $Pic(X)$ is given by $H_X,P_1,\dots,P_N$, where $H_X=\pi^{-1}L$ is the class of any line $L$ which is disjoint from all the $p_j$, and $P_j$ is the class of the divisor $\pi^{-1}p_j$.  If $C\subset {\bf P}^2$ is any curve, then we let $[C]_X$ denote its class in $Pic(X)$.  Thus $\pi^*[C]_X = m\cdot H_X +\sum \mu_j P_j$, where  $m$ denotes the degree of $C$, and $\mu_j$ is the multiplicity of $C$ at $p_j$.  (If $p_j\notin C$, then $\mu_j=0$.) 

If $f:X\to X$ is a rational map, then the pullback map $f_X^*$ is a well-defined linear map of $Pic(X)$.  We will consider $f_X^*=\left({m_{i,j}}\right)$ as a matrix with integer entries with respect to the ordered basis $H_X,P_1,\dots,P_N$.  Thus
$$f^*[L] = m_{1,1}[L] + {\rm linear\ combination\ of\ }P_1,\dots, P_N$$
\begin{proposition}  The entry $m_{1,1}$ in $f^*_X$ is the degree of $f$.
\end{proposition}

In particular, we conclude that if $(f_X^n)^* = (f_X^*)^n$, then the degree of $f^n$ is the (1,1)-entry of the matrix $(m_{i,j})^n$ and thus satisfies a linear recurrence.

Now we consider  the space $X$ obtained in the previous paragraph by blowing up $e_1$ and $e_2$.   The induced map $f^*$ on $Pic(X)$ acts according to
$$E_1\to L_\infty\to E_2\to [x+b=0]$$
Thus, $f^*:E_1\to H_X-E_1-E_2$ and $E_2\to H_X-E_2$.  

Next we need to determine what $f_X^*$ does to $H_X$.  We start by looking at ${\bf P}^2$; since $f$ has degree 2, $f^{-1}H$ is a quadric.  Both centers of blowup are indeterminate and blow up to lines.  Thus a general line $H\subset{\bf P}^2$ intersects each of these blowup images with multiplicity one, so $f^{-1}H$ is a quadric which goes through both $e_1$ and $e_2$.  In terms of divisors, this means that
$$f_X^*H_X=2H_X - E_1-E_2$$
With respect to this basis we have
$$f_X^* = \begin{pmatrix}
 2 & 1 & 1\\
-1 & -1 & 0\\
-1& -1 & -1
\end{pmatrix}$$

Let us suppose that $(a,b)\in {\cal V}_n$ and let $\pi:Y\to X$ to be the blowup of the points $q_0,\dots, q_n$ as in the previous paragraph.  Thus  $Pic(Y)=\langle H_Y,E_1,E_2,Q_n,Q_{n-1},\dots,Q_1\rangle$.  As above, the exceptional fibers are mapped as
$$f_Y: P=Q_n\to Q_{n-1}\to\cdots \to Q_1\to\{y+a=0\}$$
In terms of divisors we have $[y+a=0]_Y=H_Y-P-E_1$
and  $[x+b=0]_Y=H_Y-E_1-E_2-P$, and $f_Y^*H_Y=H_Y-E_1-E_2-P$.  The difference between $[\cdot]_X$ and $[\cdot]_Y$ arises because the curves may contain different centers of blowup.  Thus with respect to this ordered basis of $Pic(Y)$, we have
$$f_Y^* = \begin{pmatrix}
 2 & 1 & 1& && &1  \\
-1 & -1 & 0  &&&         & -1 \\
-1 & 0 & -1 & &&         &  0\\
 & & & 0 &  &      & -1\\
 & & & 1 & 0 &       \\
 &&&& 1 & 0 \\
 &&&& &  1 &0 
 \end{pmatrix}$$
\begin{proposition} The characteristic polynomial of the matrix above is 
$$\chi_n(t) = t^{n+1}(t^3-t-1) + t^3+t^2-1$$ 
If $\lambda_n$ denotes the largest root of $\chi_n$, then $\lambda_7>1$, and $\lambda_n$ is increasing in $n$.
\end{proposition}

We conclude that if $(a,b)\in{\cal V}_n$, then $\delta(f)=\lambda_n$, and thus $\delta(f)>1$ if $n\ge 7$.
\section{Matrix inversion and variations}  
Let ${\cal M}_q$ denote the space of $q\times q$ matrices, and let ${\bf P}({\cal M}_q)={\cal M}_q^*/{\bf C}^*$ denote its projectivization.  We consider the mapping $J$ defined on $q\times q$ matrices by component-wise inversion:  $J(x_{i,j})=(1/x_{i,j})$. $J$ is clearly smooth at the matrices  $x$ for which the entries are all nonzero.  We may also write $J$  as a matrix of polynomials by setting $J(x) =( x_{i,j}^{-1}\prod x)$, where $\prod x:=\prod_{(\mu,\nu)}x_{\mu,\nu}$ is the product of all of the entries of $x$.  Thus we see that $J$ has degree $q^2-1$ on ${\bf P}({\cal M}_q)$.  We let  $I(x_{i,j})=(x_{i,j})^{-1}$ be the usual matrix inversion.  Recall the familiar formula for $I(x)$ as the quotient of the classical adjoint, formed from the $(q-1)\times(q-1)$ minors, divided by the determinant.  From this we see that $I$ has degree $q-1$ as a self-map of  ${\bf P}({\cal M}_q)$.   Both of the maps  $I$ and $J$ are rational involutions, defined and regular on dense subsets of  ${\bf P}({\cal M}_q)$.  We will be concerned with the map  $K=I\circ J$ which is a birational map, and $I^{-1}\circ K\circ I=K^{-1}$, so $K$ is reversible, in the sense of being conjugate to its inverse.  
To suggest that there is subtlety in composing these maps, we note that:
\begin{proposition}  The degree of $K=I\circ J$ is $q^2-q+1 <\max({\rm deg}(I),{\rm deg}(J))$.
\end{proposition}

The map  $K$ was studied by Angl\`es d'Auriac, Maillard, and Viallet \cite{[AMV]}, as well as the restrictions of  $K$ to the subspaces ${\cal S}_q$ of symmetric matrices, and to ${\cal C}_q$ of cyclic matrices, which have the form
$$\begin{pmatrix}
a_0 & a_1 & \dots & a_{q-1}\\
&a_0 & a_1 &   \\
& \ddots &\ddots &\\
a_1& a_2 &\dots & a_0
\end{pmatrix}$$
Based on their analysis (largely numerical) of these maps, they conjectured the following:
\begin{theorem} The dynamical degrees of all three maps coincide: 
$$\delta(K)=\delta(K|_{{\cal S}_q}) = \delta(K|_{{\cal C}_q})$$
and this number is the largest root of $t^2-(q^2-4q+2)t+1$.
\end{theorem}

This Theorem was proved as a combination of results in \cite{[BT]} and \cite{[T]}.  We note that passing to a linear subspace does not increase the degree, so the inequalities $\delta(K)\ge\delta(K|_{{\cal S}_q}) $ and $\delta(K)\ge\delta(K|_{{\cal C}_q})$ follow easily.  The restriction $K|_{{\cal C}_q}$ introduces symmetries that make the map much easier to deal with. On the other hand, the additional symmetries make the restriction $K|_{{\cal S}_q}$  harder to deal with than the unrestricted $K$.  The set of symmetric, cyclic matrices ${\cal SC}_q={\cal S}_q\cap{\cal C}_q$ is also invariant under $K$.  This introduces all of the symmetries of ${\cal C}_q$ as well as ${\cal S}_q$, so there are different sorts of symmetries.  The map $q\mapsto \delta(K|_{{\cal SC}_q})$ depends on $q$ in  a more complicated way (see \cite{[BK3]}).

\section{The maps $I$, $J$ and $K$}
The maps $I$ and $J$ are involutions, so $\delta(I)=\delta(J)=1$.  We discuss the process of regularizing them by blowing up.  We define the set $\Sigma_{i,j}$ to be the set of matrices for which the $(i,j)$-entry vanishes.    Similarly, we let $e_{i,j}$ denote the matrix for which all entries are zero except in  the location $(i,j)$.  Now we consider $J$ as a map of ${\bf P}({\cal M}_q)$.  $J$ is regular at each $x=(x_{i,j})$ for which all the entries $x_{i,j}\ne0$.  We see that $J(\Sigma_{i,j}-{\cal I}(J)) = e_{i,j}$.  Conversely, since $J=J^{-1}$, we see that $J$ blows $e_{i,j}$ up to $\Sigma_{i,j}$.   Given a point $x=(x_{i,j})$, we let $T(x)$ be the set of all $(i,j)$ such that $x\in \Sigma_{i,j}$.  Then $J$ blows up $x$ to the linear subspace generated by $\{e_{i,j}: (i,j)\in T(x)\}$, which is $\bigcap _{(\mu,\nu)\notin T(x)} \Sigma_{\mu,\nu}$.  For instance, if $x_{i_1,j_1}=x_{i_2,j_2}=0$, and if all other entries of $(x_{i,j})$ are nonzero, then $J$ blows up $x$ to the line passing through $e_{i_1,j_1}$ and $e_{i_2,j_2}$.   $J$ is indeterminate at the sets  $\Sigma_{i_1,j_1}\cap \Sigma_{i_2,j_2}$ for which  $(i_1,j_1)\ne(i_2,j_2)$.  In fact,
$${\cal I}(J)=\bigcup_{(i_1,j_1)\ne(i_2,j_2)} \Sigma_{i_1,j_1}\cap \Sigma_{i_2,j_2}.     \eqno(5.1)$$

Now we define the space $\pi:X\to {\bf P}({\cal M}_q)$ in which all points $e_{i,j}\in {\bf P}({\cal M}_q)$, $1\le i,j\le q$, are blown up.   The fiber $\pi^{-1}e_{i,j}\cong{\bf P}^{q^2-2}$ is the projectivization of the normal bundle to ${\bf P}({\cal M}_q)$ at $e_{i,j}$.  (The space of tangent vectors normal to a point is the space of all tangent vectors at that point.)  That is, if $\nu$ is a vector normal to $e_{i,j}$, then the curve $t\mapsto \pi^{-1}(e_{i,j}+t\nu)$ lands at a unique point $\hat\nu\in E_{i,j}$ as $t\to0$.  The space $Pic(X)$ is spanned by the class of a general hypersurface $H_X\subset X$ and the classes of exceptional divisors $E_{i,j}$.  To define the map $J_X^*:Pic(X)\to Pic(X)$, we start with the observation that $J^{-1}E_{i,j}=\Sigma_{i,j}$, so the class $E_{i,j}$ is taken to the class of $\Sigma_{i,j}$ in $Pic(X)$.  Since the class of $\Sigma_{i,j}$ is the same as a general hypersurface $H_X$, except that it is missing the $E_{\mu,\nu}$ for all $(\mu,\nu)\ne(i,j)$, we have
$$E_{i,j}\mapsto H_X - \sum_{(\mu,\nu)\ne(i,j)}E_{\mu,\nu}.\eqno(5.2)$$
It remains to determine $J^*(H_X)$.  On ${\bf P}({\cal M}_q)$ we have $J^*H=(q^2-1)H$.  This is because if we represent $H=\sum c_{i,j} x_{i,j}$ as a linear function, then $J^*H=\sum c_{i,j} J_{i,j} = \sum_{i,j} c_{i,j} x_{i,j}^{-1}\prod x$ is represented by the linear combination of the coordinates of $J$.     At the point $e_{1,1}$, for instance, the (1,1) component of $J$ vanishes to order to $q^2-2$, and the other components vanish to order $q^2-1$.  Thus if all the $c_{i,j}$ are non-vanishing, we see that the multiplicity (order of vanishing) of $J$ at the point $e_{\mu,\nu}$ is $q^2-2$.  Thus we have
$$J^*(H_X)=(q^2-1)H_X- (q^2-2) \sum_{\mu,\nu} E_{\mu,\nu}.\eqno(5.3)$$
\begin{proposition}  (5.2--3) together determine the linear map $J_X^*$ on $Pic(X)$.  
\end{proposition}

\noindent  More details of proof can be found in \cite{[BK1]}. 

Now we discuss the map $I$ briefly.  The matrix $x={\rm diag}(0,\lambda_2,\dots,\lambda_q)\in{\bf P}({\cal M}_q)$ is mapped to $I(x)={\rm diag}(1,0,\dots,0)$.  More generally, if $x$ has rank $q-1$, then we let $v\in{\bf C}^q$ generate the kernel, and we let $w$ be an element of the dual space ${\bf C}^{q*}$  such that its kernel is the range of $x$.  It may be shown that for matrices of rank $q-1$, the inverse $I$ (projectively), interchanges kernel and range, so $I(x)=v\otimes w=(v_iw_j)$ is a matrix of rank 1.  In particular, the set $R_{q-1}:=\{x\in{\bf P}({\cal M}_q):{\rm det}(x)=0\}$ is the exceptional hypersurface for $I$, and the image $I(R_{q-1})=R_1$ is the set of matrices of rank 1.  To regularize $I$, we construct the maifold $\pi:Z\to{\bf P}({\cal M}_q)$, which blows up the set $R_1$ of rank 1 matrices.   Let ${\cal R}^1:=\pi^{-1}(R_1)$ denote the exceptional divisor.   Near the point $x_0:={\rm diag}(1,0,\dots,0)$, the set of rank 1 matrices are parametrized by $
(x_2,\dots,x_q,y_2,\dots,y_q)\mapsto \hat x^t\otimes\hat y:=(1,x_2,\dots,x_q)^t\otimes(1,y_2,\dots,y_q)$.  The fiber $\pi^{-1}x_0$ can be interpreted as the (projectivized) $(q-1)\times(q-1)$ matrices 
$\hat\xi:=\begin{pmatrix}
0 & 0& \dots & 0\\
0 & \xi_{2,2} & \dots & \xi_{2,q}\\
0 & \vdots &  & \vdots\\
0 & \xi_{q,2}& \dots & \xi_{q,q}
\end{pmatrix}$,
and a point near the fiber over $x_0$ is given by $ \hat x^t\otimes\hat y + s\hat\xi$ for some scalar $s\in{\bf C}$.
\begin{proposition} The map $I_Z:=\pi^{-1}\circ I:{\bf P}({\cal M}_q)\to Z$ is a local diffeomorphism at generic points of $R_{q-1}$.  Further, $I_Z$ is regular at all points of $R_{q-1}$ with rank  $q-1$, and $I_Z$ is a birational map from $R_{q-1}$ to ${\cal R}^1$.
\end{proposition}

Finally we turn to the map $K=I\circ J$.  Let us define ${A}_{i,j}$ to be the set of all matrices $(x_{\ell,m})$ whose entries are zero everywhere on the $i$-th row and the $j$-th column.  This is a linear subspace of $ {\bf P}({\cal M}_q)$.  We find that $K(\Sigma_{i,j})={A}_{i,j}$.  Thus we will need to work with the space $\pi:X\to {\bf P}({\cal M}_q)$ in which all the subspaces $A_{i,j}$ are blown up, and $R_1=J(R_1)$ is blown up, in addition.  We let $K_X:=\pi^{-1}\circ K\circ\pi$ be the induced map of $X$.  In the new space $X$, $\Sigma_{i,j}$ is not  exceptional for $K_X$.  Let us define the subsets ${\cal A}_{i,j}:=\pi^{-1}A_{i,j}$.  We find that $K_X$ maps ${\cal A}_{i,j}$ to $B_{j,i}:={\cal A}_{j,i}\cap \Sigma_{j,i}$.  So each ${\cal A}_{i,j}$ is exceptional.  We now construct the space $\pi:Y\to X$ in which all the subsets $B_{i,j}\subset X$ are blown up.   Working with the induced map $K_Y$ we can determine the dynamical degree $\delta(K)$.  Further details are in \cite{[BT]}.

\section{Intermediate degrees}  
In the case of projective space $X={\bf P}^k$, we let $\omega$ denote a positive, closed (1,1)-form.  Thus $\omega$ defines a K\"ahler metric on  ${\bf P}^k$.   We write the exterior powers as $\omega^p=\omega\wedge\cdots\wedge\omega$ and set $\beta_p:=\omega^p/p!\,$.  Let $M\subset {\bf P}^k$ be a compact complex submanifold of codimension $p$.   Let us normalize $\omega$ so that $\int_{{\bf P}^k}\omega^k/k!=\int_{{\bf P}^k}\beta_k = 1$.  With this normalization,  the volume of a (linear) hyperplane $H$ with respect to the metric $\omega$ is  ${\rm Vol}(H)= \int_{H} \beta_{k-1}=1$.   It is a classical result that the codimension $2p$ volume of $M$ (with respect to the metric defined by $\omega$) is given by ${\rm Vol}(M)=  \int_{M}\beta_p$.   Thus we have the identity between volume and cohomology class, and we use this to define degree in codimension $p$.  Specifically, if $L_p$ is a linear subspace of codimension $p$, then the class $\{L_p\}$ generates $H^{p,p}({\bf P}^k;{\bf Z})$, and the classes $\{L_p\} = \{\beta_p\}$ are equal.  So the class $\{M\}$ is a multiple of this class, and we use this to define the degree:  
$$\{M\} ={\rm deg}_p(M) \, \{L_p\}   {\rm \ \ where\ \ }   {\rm deg}_p(M)= \int_M \beta_p$$
This remarkable identity between degree, volume and topology serves to extend the previous definition of degree to  intermediate dimensions.

For a rational map $f:X\to Y$, there is a well-defined map on all cohomology groups $f^*:H^{p,q}(Y)\to H^{p,q}(X)$.  When $X={\bf P}^k$, we may use this to define the degree ${\rm deg}_p$ by the equation  ${\rm deg}_p(f) \, \{ \beta_p \} = f^* \{\beta_p \}$.  This is given as an integral: 
$${\rm deg}_p(f)=  \int_{{\bf P}^k} \beta_{k-p}\wedge f^*\beta_p$$
 The quantity ${\rm deg}_p$  is not invariant under conjugacy.  However, we see that ${\rm deg}_p(f\circ g)\le{\rm deg}_p(f){\rm deg}_p(g)$, so we can  define the dynamical degree as  $\delta_p(f): = \lim_{n\to\infty} \left( {\rm deg}_p(f^n) \right)^{1/n}$.  If $\varphi$ is a birational map of  ${\bf P}^k$, then we have $\delta_p(f) = \delta_p(\varphi^{-1}\circ f\circ\varphi)$.  
 
 For general $X$ is it natural to define intermediate dynamical degrees by setting
$$\delta_p(f):=\lim_{n\to\infty}||f^{n*}|_{H^{p,p}}||^{1/n}$$
In fact, if $f$ is holomorphic, then $(f^n)^*|_{H^{p,p}} = \left (f^*|_{H^{p,p}} \right)^n$.  Thus $\delta_p(f)$ is the spectral radius of $f^*|_{H^{p,p}}$.  In this case $\delta_p$ is an algebraic integer for all $p$.  It is natural to ask whether $\delta_p$ is an algebraic integer when $f$ is merely rational.   The material above was taken from  Russakovskii and Shiffman \cite{[RS]}, and the reader is invited to consult the original paper.

It is clear that the same  definition applies to meromorphic maps of complex manifolds.  In the case of a compact, K\"ahler manifold, it is classical that $p\mapsto\log \delta_p(f)$ is concave in $p$.   We have $\delta_0(f)=1$ and $\delta_k(f)\ge1$ for all maps.  Thus if $\delta_\ell(f)>1$ for some $0<\ell \le k$, the concavity implies we have $\delta_p(f)>1$ for all $0<p<k$.  

The following was obtained jointly with Jan-Li Lin:

\begin{theorem}  If $f$ is an automorphism of a compact, K\"ahler manifold, and if  $\delta_\ell(f)>1$ for some $0<\ell<k$, then $\delta_p(f)$ is irrational for all $0<p<k$.
\end{theorem}

\noindent{\it Proof. }   By the remark above, we have $\delta_p(f)>1$ for all $0<p<k$.  Let us suppose that $\delta_p(f)$ is rational.   If $f$ is an automorphism of $X$, then $\delta_p(f)$ is the spectral radius (modulus of the largest eigenvalue) of $f^*|_{H^{p,p}}$.  Since $H^{p,p}$ is an invariant subspace of $H^{2p}(X; {\bf C})$,  an eigenvalue of this restriction will also be an eigenvalue of $f^*$ acting on $H^{2p}(X;{\bf C})$.  Since $f^*$ also preserves $H^{2p}(X;{\bf Z})$ we may consider $f^*$ as a matrix with integer coefficients.  The characteristic polynomial $\chi(x)$ of $f^*$ is monic.  Thus all eigenvalues of $f^*$ are algebraic integers.  Let $\mu$ be an eigenvalue with maximum modulus.  

If $\mu$ is real, then $\mu=\pm\delta_p(f)$ is rational.  It is elementary that every rational, algebraic integer  actually belongs to ${\bf Z}$.   Now, since $f^*$ is an invertible, integer matrix, its determinant is $\pm1$.  Thus the characteristic polynomial has the form $\chi=x^m + \cdots\pm 1$.  On the other hand, since $\mu$ is an integer zero of $\chi$, $(x-\mu)$ is a factor of $\chi(x)$.  This means that  $\chi(x)= (x-\mu) p(x)= (x-\mu)(x^{m-1} + \cdots+ c_0)= x^m+\cdots-\mu c_0 = x^m+\cdots\pm 1$.  This is not possible since $c_0$ is an integer, and $|\mu|>1$.

If $\mu$ is not real, then we have $|\mu| = |\mu\bar\mu|^{1/2}=\delta_p(f)$, which is assumed to be rational.  Now let $\alpha_3,\dots,\alpha_m$ denote the other roots of $\chi$.  Since these are algebraic integers, it is elementary (see \cite{[M]}) that their product $\alpha_3\cdots\alpha_m$ is also an algebraic integer.  Since $\mu\bar\mu\alpha_3\cdots\alpha_m=\pm1$, we conclude that both $\mu\bar\mu$ and $\alpha_3\cdots\alpha_m$ are rational.  Since, in addition, these are both algebraic integers, they both are integers.  But this contradicts the assumption that $|\mu|>1$.

\section{Monomial maps} 
The intermediate dynamical degrees are important for understanding the dynamical behavior.   They are invariant under birational conjugacies in the following strong sense:  If $\varphi:X\to Y$ is birational, and if $g:=\varphi^{-1}\circ f\circ\varphi$, then  $\delta_p(f,X)=\delta_p(g,Y)$ (see \cite{[DS]}).  In the same paper, Dinh and Sibony give an estimate on the topological entropy of $f$:
$$h_{\rm top}(f)\le \log\max(\delta_1(f),\dots,\delta_k(f) )$$
In case $f$ is holomorphic, this is known to be an equality.   And if $f$ is holomorphic, then $f^*$ on $H^{p,p}$, is represented by an integer matrix.  The degree $\delta_p$ will be the spectral radius of this matrix and thus an algebraic integer.  On the other hand,   it is a different matter to try to find $\delta_p$ for maps which do not satisfy $(f^*)^n=(f^n)^*$ on $H^{p,p}$.

So far, the only nontrivial class on which $\delta_p$ has been computed is the monomial maps.  Let $A=(a_{i,j})$ be a $k\times k$ matrix with integer entries.  We let 
$$f_A(x)=\left (\prod_j x_j^{a_{1,j}} , \dots, \prod_j x_j^{a_{n,j}} \right)$$ 
be the monomial map defined by $A$.  It is easily seen that $f_A^n = f_{A^n}$, so the iterates are easily given.  Further, $f_A$ is a well defined rational map of ${\bf P}^k$,  and $f^*_A[L_p]={\rm deg}_p(f_A)[L_p]$.  In fact, this number is given by an integral: ${\rm deg}_p(f) =   \int \beta_{k-p}\wedge f^*\beta_p$.  The number $\delta_p$ would then be the limit of $({\rm deg}_p(f^n))^{1/n}$ as $n\to\infty$.  Although this approach is simple to describe, it seems not to be so simple to carry out.

A useful approach to finding the number $\delta_p$ in the case of monomial maps is to change the space $X={\bf P}^k$ to the space $Y=({\bf P}^1)^k={\bf P}^1\times\cdots\times{\bf P}^1$, which is birationally equivalent to  $X$.  We may let $[x_j:y_j]$ be homogeneous coordinates on the $j$-th factor of ${\bf P}^1$.  Then a basis for $H^{p,p}$ is given by the classes $L_I=\{x_{i_1}=\cdots=x_{i_p}=0\}$, where $I=(i_1,\dots,i_p)$ is a $p$-tuple of indices $1\le i_j<\dots<i_p\le k$.  (Of course, these are the same as the classes $\{\zeta_{i_1}=\cdots=\zeta_{i_p}=0\}$, where each $\zeta_j$ is either $x_j$ or $y_j$.)  We consider $\{L_I\}$ as an ordered basis for $H^{p,p}(Y)$.   
Given a matrix $M=(m_{i,j})$ let us use the notation $|M|:= \left(|m_{i,j}| \right)$ for the matrix consisting of the absolute values of the entries of $M$.  The action of $f_A^*$ on $H^{p,p}(Y)$ now has a simple description (see \cite{[L]}):
\begin{proposition} Let $M:=\bigwedge^pA$ denote the $p$-th exterior power of the matrix $A$.  Then when we write the basis $\langle L_I\rangle$ suitably, the action $f^*_A|_{H^{p,p}}$ is given by $|M|$. 
\end{proposition}

While we are working with $({\bf P}^1)^k$, it is useful to consider the degree as the matrix ${\rm Deg}_p(f)$ which represents $f^*_{H^{p,p}}$.  For instance,  $A=\begin{pmatrix} 1 & -1\\ -2&-3\end{pmatrix}$, so we have 
$f_A(x_1,x_2) = (x_1/x_2, x_1^{-2}x_2^{-3})$.  In homogeneous coordinates, this becomes
$$f_A: \ \ [x_0:x_1:x_2]\mapsto [x_1^2 x_2^3: x_1^{3}x_2^2:x_0^5]$$
so ${\rm deg}_1(f_A) = 5$, and ${\rm Deg}_1(f_A) = \begin{pmatrix}1&1\\ 2&3\end{pmatrix}$.

Now let us write the eigenvalues of $A$ as $\mu_1,\dots,\mu_k$, where $|\mu_1|\ge|\mu_2|\ge\cdots\ge|\mu_k|$.  The following result, obtained independently by C. Favre and E. Wulcan \cite{[FW]},  and J-L Lin \cite{[L]}, gives the dynamical degrees:
\begin{theorem}  The dynamical degrees are $\delta_p(f_A) = |\mu_1\cdots\mu_p|$, $1\le p\le k$.
\end{theorem}

The idea of why the Theorem follows from the Proposition is as follows.  The exterior product is  $(\bigwedge^p A)(v_1\wedge\cdots\wedge v_p):=(Av_1)\wedge\cdots\wedge (Av_p) $.  If  $v_i$ is an eigenvector satisfying $Av_i = \mu_i v_i$, then $(\bigwedge^pA)(v_1\wedge \cdots\wedge v_p) = (\mu_1\cdots\mu_p)v_1\wedge \cdots\wedge v_p$.  The size of $\bigwedge ^p (A^n)$, and thus $\left|\bigwedge^p(A^n)\right|$, can be estimated above and below by $|\mu_1\cdots\mu_p|^n$, which gives the claimed exponential growth.

\end{document}